\documentclass[a4paper]{article}

\usepackage{amsmath,amstext}  
\usepackage{amsfonts,amssymb, bbm} 
\usepackage[LGR, T1]{fontenc}
\usepackage[utf8]{inputenc}
\usepackage[english]{babel}
\usepackage{color}
\usepackage{graphicx, graphics}
\usepackage{subfigure}
\usepackage{epsfig}
\usepackage{epsf}
\usepackage{vmargin}
\usepackage{enumerate}
\usepackage{tikz}
\usepackage{dsfont}
\usepackage{bbm}
\usepackage{stmaryrd}
\usepackage{multirow}
\usepackage{algorithm,algorithmic}

\begin{document}


\title{Optimization of a launcher integration process: a Markov decision process approach\thanks{This work was supported by Airbus DS and the Conseil R\'{e}gional d'Aquitaine.}}

\author{Christophe Nivot\thanks{Inria Bordeaux Sud-Ouest, Universit\'{e} de Bordeaux, IMB, CNRS UMR 5251, France} \and Beno\^{\i}te de Saporta\thanks{Universit\'{e} de Montpellier, IMAG, CNRS UMR 5149, France} \and Fran\c{c}ois Dufour\thanks{Bordeaux INP, IMB, CNRS UMR 5251, Inria Bordeaux Sud-Ouest, France} \and Damien B\'{e}rard-Bergery\thanks{Airbus DS} \and Charles Elegbede\thanks{Airbus DS}}

\maketitle

\begin{abstract}
This paper is dedicated to the numerical study of the optimization of an industrial launcher integration process. It is an original case of  inventory-production system where a calendar plays a crucial role. The process is modeled using the Markov Decision Processes (MDPs) framework. Classical optimization procedures for MDPs cannot be used because of specificities of the transition law and cost function. Two simulation-based algorithms are tuned to fit this special case. We obtain a non trivial optimal policy that can be applied in practice and significantly outperforms reference policies.
\end{abstract}
\section{Introduction}
The general class of inventory-production systems is often associated to cost optimization problems. Indeed, one must deal with three major matters: the storage of components, the possible random behavior of the manufacturing process and random clients' demand \cite{JM74,CL88,SL96}. The controller must decide which production rate of the components fits best. A too slow production rate leads to low stock levels but it might not meet clients' demand. On the opposite, a fast production rate does meet the demand, but may raise stock levels. One must then find a balance between both to minimize costs. 

This paper focuses on the optimization of a real-life industrial launcher integration process studied in collaboration with Airbus Defence and Space. Clients order a certain number of launches to be performed at specific dates. The controller has to determine the production rates in order to minimize costs. Only storage and lateness costs are taken into account here. In general, the costs may also take into account several other constraints such as exploitation cost, workforce salary, the cost related to the unavailability of the structure including any penalty or the maintenance and inspection cost, among others. Plus, a part of the architecture of the process is not set. Indeed, the controller has to decide on the maximum capacity of one warehouse between two options. 

The originality of this problem is twofold. On the one hand, the optimization horizon is rather long, 30 years, but the controller can only make decisions once a year concerning the production rates. On the other hand, the launches must be performed according to a prescribed calendar corresponding to clients' orders.

Our goal is to find an optimization procedure usable in practice. It should provide explicit decision rules applicable to each trajectory as a table giving the controller the best action to take according to the current state and time. A preliminary study was performed on a simplified process \cite{E13} using Petri nets \cite{VDA94,J97}. Although they are easy to simulate, they are not suitable for performing dynamic decisional optimization. A more suitable framework is that of Markov Decision Processes (MDPs) \cite{P94,HLL96,BR11}. MDPs are a class of stochastic processes suitable for cost and decision optimization. Briefly, at each state, a controller makes a decision which has an influence on the transition law to the next state and on a cost function. The latter depends on the starting state and the decision made. The sequence of decisions is called a policy, and its quality is gauged thanks to a cost criterion (typically, it is the sum of all the costs generated by the transitions). 

The first step to solve our problem is to implement an MDP-based simulator of the launcher integration process. Simulation results were presented at the ESREL conference in 2015 \cite{N15}. This paper deals with the optimization itself. It is a non standard optimization problem within the MDP framework because the transition law is not analytically explicit, it is only simulatable. Thus, standard optimization techniques for MDPs such as dynamic programming \cite{H60,B87,P94}, or linear programming \cite{P94,HLL96} do not apply. In addition, the cost function is unusual as the actual lateness can be computed only at the end of a year, and not at its beginning when the controller makes their decisions.

As the launcher integration process can be simulated, we investigate simulation-based algorithms for MDPs \cite{HS04,CFHM07,CHFM13}. These extensively use Monte-Carlo methods to estimate the performance of a policy. Thus, they require a fast enough simulator for the algorithms to give a result within a reasonable time. New difficulties arise here. First, the state space of our MDP, though finite is huge. Second, the first simulator in MATLAB is not fast enough. Third, the algorithms require the computation of a product of numerous numbers between $0$ and $1$, and although the output is non zero on paper, it is treated as zero numerically, leading to erroneous results. To overcome these difficulties, we reduce the state space by aggregating states in a manner that makes sense regarding our application, we use the C language and a special logarithmic representation of numbers. The results we obtained are presented and discussed.

This paper is organized as follows. Section \ref{lauint} is dedicated to the description of the assembly line under study and the statement of the optimization problem. In section \ref{mardec}, we present how the optimization problem for the assembly line fits into the MDP framework. Section \ref{optlau} presents the main difficulties encountered while trying to optimize our MDP, and solutions to bypass them. In section \ref{numres}, we present and comment the numerical results obtained. Finally a last section gives some concluding remarks. Technical details regarding the implementation of algorithms are provided in the Appendix.
%
\section{Launcher integration process description}\label{lauint}
\begin{figure*}
\begin{center}
\includegraphics[scale=0.35]{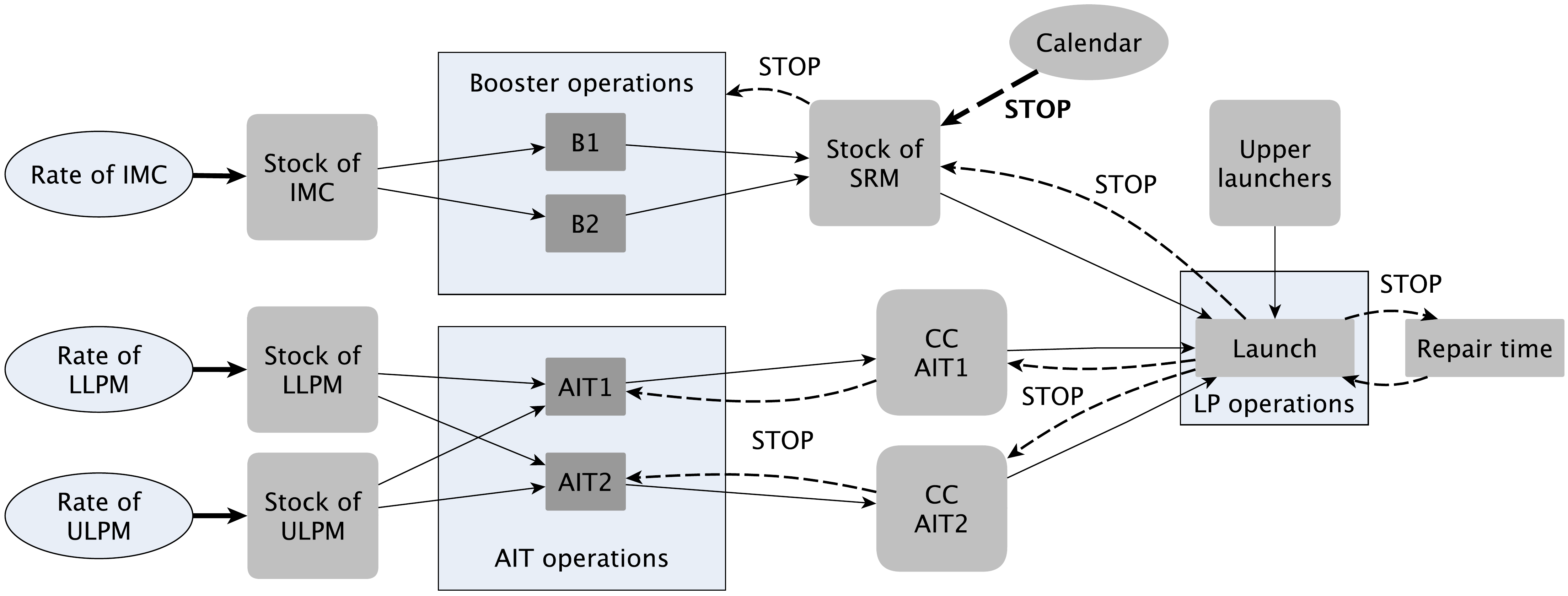}
\end{center}
\caption{The launcher integration process}
\label{Process}
\end{figure*}
Airbus Defense and Space (Airbus DS) as prime contractor is in charge of launchers and ground facilities design.
This paper is dedicated to the optimization of an assembly line representative of a launcher integration process managed by Airbus DS.
For confidentiality matters, all parameter values and random distributions given in this paper are arbitrary but realistic.

The launcher integration process we study in this paper is depicted on figure \ref{Process}. 
This assembly line is composed of several workshops and storage facilities that are described in detail in the following sections, and is operational typically for 30 years.
%
\subsection{Subassemblies}
The subassemblies are the input of the assembly line.
A launcher needs four types of subassemblies to be manufactured. These are
\begin{itemize}
\item the Insulated Motor Cases (IMCs), which are powder-free boosters,
\item the Lower Liquid Propulsion Modules (LLPMs) and the Upper Liquid Propulsion Modules (ULPM), which are the lower part of a launcher,
\item the Upper Launchers, which are the fairings of the launchers.
\end{itemize}

The Upper Launchers are always available when needed. So in the following, their production rate and storage cost are not taken into account.

The production time of the other subassemblies (IMC, LLPM, ULPM) is assumed to be random. 
Its distribution is symmetric around a mean $\tau$ and takes 5 different values with probabilities given in Table \ref{lawsub}.
The average production time $\tau$ is computed by dividing the total number of workdays by the target number of subassemblies the controller decides to produce in a year, and taking its integer part.
The number of workdays in a year is set at 261. So for instance, if the controller wants an average production of 12 LLPMs a year, the value of $\tau$ for LLPMs will be $[261/12] = 21 $ days, where $[x]$ stands for the integer part of a real number $x$.
\begin{table}[tp]
\centering
\begin{tabular}{l l}\hline
Value (days) & Probability\\
\hline
$\tau-2$ & 3/32\\
$\tau-1$ & 5/32\\
$\tau$ & 1/2\\
$\tau+1$ & 5/32\\
$\tau+2$ & 3/32\\
\hline
\end{tabular}
\label{lawsub}
\caption{Distribution of the production time of a subassembly with mean $\tau$}
\end{table}

Each produced unit is transferred to its dedicated warehouse. Each warehouse capacity is limited to 4 units. When full, the production of the corresponding subassembly stops. It resumes as soon as one unit is taken from the stock. The subassemblies leave their warehouse when they are needed in workshops. Otherwise, they wait there.

Storage is costly. For each stored subassembly, its daily amount is a percentage of the price of the subassembly. All storage prices are presented in Table \ref{storage} as percentages of a certain reference value.
\begin{table}[tp]
\centering
\begin{tabular}{l l}\hline
Subassembly & Storage price apiece per day\\
\hline
IMC & 2.6\%\\
LLPM & 55.94\%\\
ULPM & 35.59\%\\
SRM & 8.08\%\\
CC & 100\%\\
\hline
\end{tabular}
\label{storage}
\caption{Storage price for each subassembly}
\end{table}
%
\subsection{Workshops}
The process comprises various types of operations corresponding to dedicated workshops:
\begin{itemize}
\item the Booster operations refer to the integration of the boosters,
\item the Assembly, Integration and Test (AIT) operations refer to the integration of the lower part of the launcher,
\item the Launch Pad (LP) operations refer to the integration of the final launcher and the launch.
\end{itemize}
Workshops for Booster and AIT operations comprise two docks working in parallel. So they can take two subassemblies at the same time for integration. 
The Launch Pad has a single dock. 
The workshops are subject to breakdowns, maintenance operations or staff issues. Thus, their operating time is random with a uniform distribution on small set of values around the nominal values given in Table \ref{lawopt}.
\begin{table}[tp]
\centering
\begin{tabular}{l l}\hline
Workshop & Set of possible values\\
\hline
Booster docks & 5, 5.5\\
AIT docks & 25, 25.5, 26\\
Launch Pad & 10, 10.5\\
\hline
\end{tabular}
\label{lawopt}
\caption{Law of the operating time for each workshop}
\end{table}

The Booster operations start with the production of one IMC. It goes to an available dock if any. Nominal operating time for the integration of IMCs is 5 days. When integration is over, the output is called a Solid Rocket Motor (SRM). The SRM is brought to its dedicated warehouse which has a limited capacity of 4 or 8 units. The storage cost is given in Table \ref{storage}. It is up to the controller to decide which capacity fits best. When the SRM stock is full, docks B1 and B2 (Fig. \ref{Process}) can no longer work and must wait until SRMs are taken away from it. SRMs leave their warehouse to go to the Launch Pad by groups of 4 units.

The AIT operations require one LLPM and one ULPM to be performed. When the resources are sufficient, they go to an available dock if any. Nominal operating time for their integration is 25 days. The output of these operations is called a Central Core (CC). When a CC is manufactured, it waits inside its dock until it is needed in the LP. During this time, the dock is occupied and cannot be used for the integration of another CC. Keeping a CC in its dock is treated like storage, hence is costly, see Table \ref{storage}.

The LP operations need one Upper Launcher (which is always available by assumption), one CC and 4 SRMs. Nominal operating time for the integration and the launch is 10 days. The specificity of the Launch Pad is that after a launch is performed, it goes through a series of repair operations which invariably last 5 days. It is the only duration which is not random for this process. During the repair time, the LP is not available and next subassemblies have to wait until the end of these five days to be integrated there.
%
\subsection{Calendar}
A major specificity of this assembly line is the launch calendar. The launches have to be performed according to a predetermined schedule requested by clients. As a first approach, we make the simplifying assumption that it is established once and for all for 30 years. In real life operations, the launch calendar is usually known two years in advance. An admissible calendar has to respect the following constraints.
On the one hand, the first four years are dedicated to the system's starting. During this time, 18 launches have to be performed:
\begin{itemize}
\item year 1: 1 launch,
\item year 2: 2 launches,
\item year 3: 4 launches,
\item year 4: 11 launches.
\end{itemize}
On the other hand, for each of the next 26 years, the number of launches is between 6 and 12. The launches are closer to each other at the end of a year than at its beginning. In average, the number of launches per year is around 10. For security reasons, launches have to be spaced out in time by at least 15 days.

The calendar has a direct effect on the integration process itself. Indeed, the launch must not be ahead of schedule. If the next launch date is too far away in time (typically 10 days away or more), the integration process stops. It can be modeled by blocking the stock of SRM, therefore delaying the next launch. It gets unblocked as soon as the launch date is close enough.

\begin{table}[tp]
\centering
\begin{tabular}{l l}\hline
Type of lateness & Price per day\\
\hline
Anticipated & 45.19\%\\
Unexpected & 80.13\%\\
\hline
\end{tabular}
\label{lateness}
\caption{Lateness costs}
\end{table}
However, a launch may occur later than originally scheduled. Lateness is costly. There are two types of lateness costs. Their price are proportional to their duration and are reported in table \ref{lateness}.
\begin{itemize}
\item When the Launch Pad starts to integrate a launcher less than 10 days before the launch date, lateness is certain but anticipated.
\item If the Launch Pad started on time, but the launch is late all the same, it is called an unexpected lateness.
\end{itemize}
Note that it is much more costly to be late without anticipation.
%
\subsection{Objective}
This assembly line is monitored by a controller. There are two types of decision to be made.
\begin{itemize}
\item The capacity of the SRM warehouse must be fixed for the whole process (static decision). There are two options: 4 or 8 units. A smaller capacity costs less but might slow down the process. The decision is made before the process starts. 
\item The production rate of IMCs, LLPMs and ULPMs are decided at the beginning of each year, essentially for logistic reasons (dynamic decision). The production rate is the average number of subassembly produced in a year. For the IMCs, it goes from 24 to 48 by a step of 4. For the LLPMs and the ULPMs, it goes from 6 to 12. Slowing down their production may be useful to keep the storage levels reasonable, so that they do not cost too much. However, a low production rate will generate lateness.
\end{itemize}
The cost of the process is computed at the end of every year, and all these annual costs are added up to make the total cost. The objective of this paper is to determine the best policy for choosing the rates of production and the architecture to meet the demand and minimize the total cost. This policy should also be numerically tractable and applicable in practice.
%
\section{Markov decision processes model}\label{mardec}
This section contains a brief presentation of the basics of Markov decision processes (MDPs)\cite{HLL96} and how they can be used to solve our optimization problem for the launcher integration process.
\subsection{Discrete-time Markov decision processes}
\begin{figure}
\begin{center}
\begin{tikzpicture}
\node[draw,rounded corners=3pt,text width=3.6em,fill=gray!20] (etatn) at (0,0) {$x_t$ state at time $t$};
\node[draw,rounded corners=3pt,text width=3.7em] (agent) at (1.8,1) {$a_t$ decision at time $t$};
\node[draw,rounded corners=3pt,text width=3.6em] (cout) at (2.3,-1) {Cost $c(x_t,a_t)$};
\node[text width=5.8em] (trans) at (4.5,0) {Transition probability $Q(\cdot \mid x_t,a_t)$};
\node[draw,rounded corners=3pt,text width=4.2em,fill=gray!20] (etatn1) at (6.7,0) {$x_{t+1}$ state at time $t+1$};

\draw[->,>=latex] (1.8,0.38) -- (1.8,0);
\draw[->,>=latex] (2.3,0) -- (2.3,-0.5);
\draw[-,>=latex] (etatn) -- (2.5,0);
\draw[-,>=latex,dashed] (2.5,0) -- (trans);
\draw[->,>=latex,dashed] (5.1,0) -- (etatn1);
\end{tikzpicture}
\end{center}
\caption{Principle of a Markov decision process}
\label{MDP}
\end{figure}
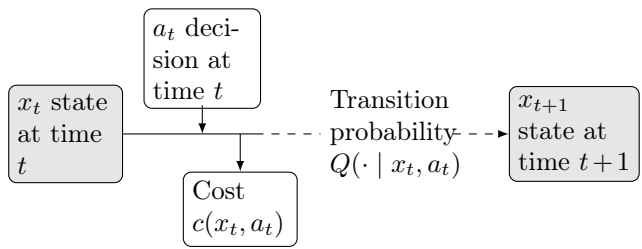

MDPs are a class of stochastic controlled processes. They model processes subject to randomness which trajectory can be influenced by external decisions. A standard MDP is presented in figure \ref{MDP}.

A discrete-time MDP is a sequence of states $(x_t)$ called a trajectory. It represents the succession of the states the system is in. 
At each step $t$, the system is in state $x_t$ belonging to a set called state space. It contains all the possible states. A decision $a_t$ is made according to the whole history of the system (states and decisions). It belongs to a set called action space. The admissible decisions may depend on the state $x_t$. Once an action is decided, a cost depending on $x_t$ and $a_t$ is generated. It is called a cost-per-stage. Then, the system goes to state $x_{t+1}$ according to a probability transition depending on $x_t$ and $a_t$. That is why this class of processes is said to be controlled.

Optimization problems associated to MDPs use cost criteria to determine which decisions are the best. These criteria consist in adding up all the costs-per-stage. The aim is to minimize the cost criterion over the admissible sequences of actions, called policies. Policies are said to be acceptable if they lead to perform all the prescribed launches within the 30 years.
%
\subsection{Modeling the assembly line as an MDP}
The complexity of the launcher integration process raises problems of time scaling. The natural time scale to describe the process and compute the costs is half-days. The natural time step to take decisions is 1 year. Thus, the model includes these two time scales.

To simulate the launcher integration process with accuracy yet another discrete time scale is required. It is called the event scale. The transition from one state to another corresponds to a major event such as the production of a subassembly or a launch for instance. Possible events are:
\begin{itemize}
\item the production of an IMC,
\item the production of an LLPM,
\item the production of a ULPM,
\item either dock B1 or dock B2 produced an SRM,
\item either dock AIT1 or dock AIT2 produced a CC,
\item the LP performed a launch or finished repairing,
\item the calendar authorizes the SRM storage to be unblocked,
\item the current year has ended.
\end{itemize}
To determine which one of these events comes next, one integrates additional components in the state variable leading to a 28-dimensional state variable. For instance, we added all the remaining times before each event happens. Finding their minimum leads to selecting the upcoming event. These variables are not observable in practice, they are useful for simulation purposes only. 
To adapt this event time scale to the yearly decision-making time scale, we introduced the event "the current year has ended". This way, one can obtain the state variable at the beginning of a year. The succession of such yearly states is the MDP we are working with. This is the right time scale for optimization.
%
\section{Optimization of the launcher integration process}\label{optlau}
%
There are many optimization techniques to solve problems involving MDPs. These are, for instance, dynamic programming, linear programming, but also policy iteration \cite{H60,B87,WW89,D03} or value iteration \cite{B87,HLL96}. However, the MDP we obtain is non standard regarding its optimization. Despite being simulatable, one cannot write its transition law explicitly. Plus, the cost-per-stage function does not have the usual form. Indeed, one knows the cost generated in a year only once the year is over. Thus, one cannot use standard techniques for MDP optimization cited above because they are founded on the complete knowledge of these two elements. 

Instead, we investigated simulation-based algorithms \cite{CFHM07,CHFM13,HFM08,HHC12}. Yet, there is another difficulty: the state space is huge. One may estimate the total number of states around 9 billion. Thus, even simulation-based algorithms may be computationally inefficient. This is why a first step to reduce the state space is needed.
%
\subsection{Reduction of the state space}
As mentioned above, the controller does not have access to all the variables used for the simulation of the launcher integration process. Not all of them are directly observable and among them, only a few are relevant to the decision making. These are the stock levels and the number of launches to be performed in the year. The state space is then reduced to six state variables. Its number of elements is then:
\begin{itemize}
\item 26,251 if the SRM storage capacity is 4,
\item 47,251 if the SRM storage capacity is 8.
\end{itemize}

For memory size issues, we also aggregated some states:
\begin{itemize}
\item In a given year, the number of launches to be performed can be greater than 12 because of late launches from the previous year. It is gathered in one state called "12 launches and more".
\item The stock level of IMCs, LLPMs, ULPMs and SRMs is represented by an integer from 1 to 3 (empty/not enough for one launcher/enough for one launcher/full), see Tables \ref{codeILU} and \ref{codeSRM}.
\end{itemize}
\begin{table}[tp]
\centering
\begin{tabular}{l l}\hline
Number of units in store & Representation\\
\hline
0 & 1\\
1, 2 or 3 & 2\\
4 & 3\\
\hline
\end{tabular}
\label{codeILU}
\caption{Aggregation of the IMCs, LLPMs and ULPMs stock levels}
\end{table}
\begin{table}[tp]
\centering
\begin{tabular}{l l}\hline
Number of units in store & Representation\\
\hline
0, 1, 2 or 3 & 1\\
4, 5, 6 or 7 & 2\\
8 & 3\\
\hline
\end{tabular}
\label{codeSRM}
\caption{Aggregation of the SRMs stock levels}
\end{table}
Note that we tested several representations for the stocks. This one led to the best results.

Thus, the total number of states falls down to 2,281, and the number of possible policies to ${343}^{{68,430}}$. It is still huge, which makes this optimization problem still not trivial, even if it is now numerically tractable.
%
\subsection{Simulation-based algorithms}
We investigated two simulation-based algorithms \cite{CFHM07,CHFM13}: Model Reference Adaptive Search (MRAS) and Approximate Stochastic Annealing (ASA). This section is dedicated to their informal description, see Algorithms~\ref{MRASa} and \ref{ASAa} in Supplementary Material for details. Basically, such algorithm simulate policies randomly and try to select the best ones.

Both algorithms operate on a 3-dimensional matrix $P$. Each coefficient $P(i,j,t)$ stands for the probability to take the decision $a_j$ at time $t$ when the system is in state $x_i$. 
At each iteration, the algorithms simulate trajectories, compute their cost and update the matrix $P$ accordingly. 
The MRAS algorithm updates the matrix $P$ according to a performance threshold, only when it is above the threshold. The ASA algorithm computes a weighted average of all simulated trajectories and therefore updates $P$ ate every step.
At the end, the algorithms returns an optimal matrix $P^*$ leading to a generator of the optimal policy. Indeed, a policy generated with $P^*$ is a two-dimensional matrix $M^* = (m_{i,t})$ with values in $\{1,2,\ldots,343\}$ where each $m_{i,t}$, $1 \leq i \leq 2281$, $1 \leq t \leq 30$, is the number of the action to be taken when in state number $i$ on year number $t$ according to this policy. Note that these policies are naturally dynamic and path-adapted. They can be easily applied in practice as the controller just needs to look up in the matrix $M^*$ to select the next decision optimally.

Implementing these algorithms raises two more difficulties. On the one hand, the components of matrix $P$ are probabilities, hence numbers between $0$ and $1$. They are computed
as the product of typically 23,471,490 numbers between 0 and 1 (2,281 states by 343 decisions by 30 years). This product should be positive but is below the machine accuracy. Thus, it is numerically equal to zero. This causes numerous errors in the output of the algorithms. To fix this problem one uses a specific representation of numbers using a logarithm transform.

On the other hand, the algorithm may lose some precious time computing policies leading to enormous cumulated lateness on launches, for instance if selecting the lowest rate of 6 launchers per year when 12 launches are planned in the calendar. In order to avoid these policies as much as possible, a high penalty cost of 10,000,000 is set for every launch that was not performed at the end of the 30 years.
%
\section{Numerical results}\label{numres}
The results of our numerical computations are now presented. The computer used is a Mac Pro 2.7 GHz 12-Core with 64 Go RAM.

First, we elaborated an MDP-based simulator using MATLAB. 
Preliminary tests of the MRAS and the ASA algorithms required to simulate about 75,000,000 trajectories. Our MATLAB simulator simulates one trajectory in 0.7 second. This is too long for the tests. So we implemented our simulator in C language. The C simulator simulates a trajectory in 0.003 second. Combined with parallel computing, our tests could be done in tractable time.
%
\subsection{Na\"{\i}ve policies}
%
As the optimisation problem we consider does not have a closed mathematical solution, one needs some criteria to estimate the degree of optimality of output policies from the MRAS and ASA algorithms. The natural way to compare policies is to compare their cost, the one with minimal cost being the best. 

As a reference, we introduce the so-called na\"{\i}ve policies. 
In a given year, they prescribe to manufacture the exact number of subassemblies (on average) required for the number of launches planned in the calendar for this year. For instance, if $10$ launches are planned a given year, the na\"{\i}ve controller chooses to produce on average 40 IMC, 10 LLPM and 10 ULPM. This strategy is deemed na\"{\i}ve as it does not take into account the available stocks of subassemblies at the beginning of the year. If such stocks are high, one may manage to perform the prescribed launches with slower rates. Thus, one expects the MRAS and ASA algorithms to compute policies that can outperform the na\"{\i}ve one.
%
\subsection{Comparison of the MRAS and ASA algorithms}
The aim of the first tests is to tune the parameters of the MRAS ans ASA algorithms and check that they can provide policies that outperform the na\"{\i}ve one.
For these tests, we assume that:
\begin{itemize}
\item the horizon time is 10 years,
\item the SRM storage capacity is 8 units,
\item the calendar demands 10 launches from year 5 to year 10, uniformly distributed over the year,
\item the rate of production of the subassemblies is reduced to between 8 and 12 launchers per year.
\end{itemize}

Table \ref{naive} recapitulates the estimation of the performance of the policies requiring the same fixed production rates over the years 5 to 10. The na\"{\i}ve policy described above corresponds to 40 IMC, 10 LLPM and 10 ULPM. 
As one may have expected, it is the least costly policy among the other fixed-rate ones. Indeed, manufacturing less subassemblies generates lateness, producing more of them generate useless storage costs. 
\begin{table}[tp]
\centering
\begin{tabular}{l l l l}\hline
\multicolumn{3}{c}{Average number of units per year} & \multirow{2}{*}{Performance}\\
IMC & LLPM & ULPM & \\
\hline
32 & 8 & 8 & 123,770,000\\
36 & 9 & 9 & 45,666,000\\
40 & 10 & 10 & 809,540\\
44 & 11 & 11 & 945,340\\
48 & 12 & 12 & 972,440\\
\hline
\end{tabular}
\label{naive}
\caption{Estimation of the performance of na\"{i}ve policies computed with a Monte Carlo method with 100,000 simulations}
\end{table}

The tests presented here with the MRAS were performed with the following parameters (see algorithm \ref{MRASa} in Supplemental Material): the initial number of policies to simulate is $N_0 = 100$, the initial number of simulations for the Monte Carlo method $M_0 = 1000$, $\mu = 10^{-8}$, the initial quantile $\rho_0 = 0.25$, $\alpha = 1.02$, $\beta = 1.0205$, $\lambda = 0.4$, $\nu = 0.5$ and $\varepsilon = 1$. We performed tests with other values of these parameters, these ones yield the best results we obtained.
The parameters used for the tests with the ASA are the following (see algorithm \ref{ASAa} in Supplemental Material): the initial number of policies to simulate is $N_0 = 100$, the initial number of simulations $M_0 = 5000$, the initial temperature $T_0 = 2$.

As one can see on table \ref{perfcomp}, the na\"{i}ve policy is not optimal. Indeed, the ASA has found a better one. Moreover, it has found it in less time than the MRAS did: the MRAS needed 11 hours to give these results, whereas the ASA needed 2 hours. Thus, the ASA algorithm seems to be more reliable and will be chosen for further investigations. The optimal policies computed by the two algorithms do not present particularities that could make them easy to explain. 
\begin{table}[tp]
\centering
\begin{tabular}{c c c}\hline
& $K = 100$ & $K = 150$\\
\hline
MRAS & 797,656 & 780,054\\
ASA & 727,136 & 724,899\\\hline
\end{tabular}
\label{perfcomp}
\caption{Performance of the best policy found by the algorithms with respect to the number $K$ of iterations}
\end{table}
%
\subsection{30-year optimization}
Turning back to the original problem, one wants to find the best rates of production given a 30-year calendar. One also wants to select the best SRM storage capacity. For comparison's sake, we will also give the performance of the corresponding na\"{i}ve policy.

To generate calendars randomly, the number of launches in a year is drawn according to the distribution described in Table \ref{lawlaunches}. Its expectation is 9.77 and its standard deviation is 2.80. Given the number of launches, launch dates are fixed within the year and spread out as in Table \ref{calendars}.
\begin{table}[tp]
\centering
\begin{tabular}{l l}\hline
Number of launches per year & Probability\\
\hline
6 & 1/16\\
7 & 1/16\\
8 & 1/12\\
9 & 3/24\\
10 & 1/3\\
11 & 1/6\\
12 & 1/6\\
\hline
\end{tabular}
\label{lawlaunches}
\caption{Law of the number of launches per year}
\end{table}
\begin{table*}
\centering
\begin{tabular}{l l l l l l l l l l l}\hline
Number of launches & 1 & 2 & 4 & 6 & 7 & 8 & 9 & 10 & 11 & 12\\
\hline
Launch 1 & 130 & 87 & 52 & 37 & 32 & 29 & 27 & 26 & 23 & 21\\
Launch 2 & & 174 & 104 & 74 & 64 & 58 & 54 & 52 & 46 & 42\\
Launch 3 & & & 156 & 111 & 96 & 87 & 81 & 78 & 69 & 63\\
Launch 4 & & & 208 & 148 & 128 & 116 & 111 & 107 & 92 & 84\\
Launch 5 & & & & 185 & 160 & 145 & 136 & 129 & 121 & 117\\
Launch 6 & & & & 222 & 192 & 174 & 161 & 151 & 141 & 135\\
Launch 7 & & & & & 224 & 203 & 186 & 173 & 161 & 153\\
Launch 8 & & & & & & 232 & 211 & 195 & 181 & 171\\
Launch 9 & & & & & & & 236 & 217 & 201 & 189\\
Launch 10 & & & & & & & & 239 & 221 & 207\\
Launch 11 & & & & & & & & & 241 & 225\\
Launch 12 & & & & & & & & & & 243\\
\hline
\end{tabular}
\label{calendars}
\caption{Launch dates according to the number of launches to be performed}
\end{table*}
\subsubsection{First calendar}
\begin{figure}[tp]
\begin{center}
\includegraphics[scale=0.45]{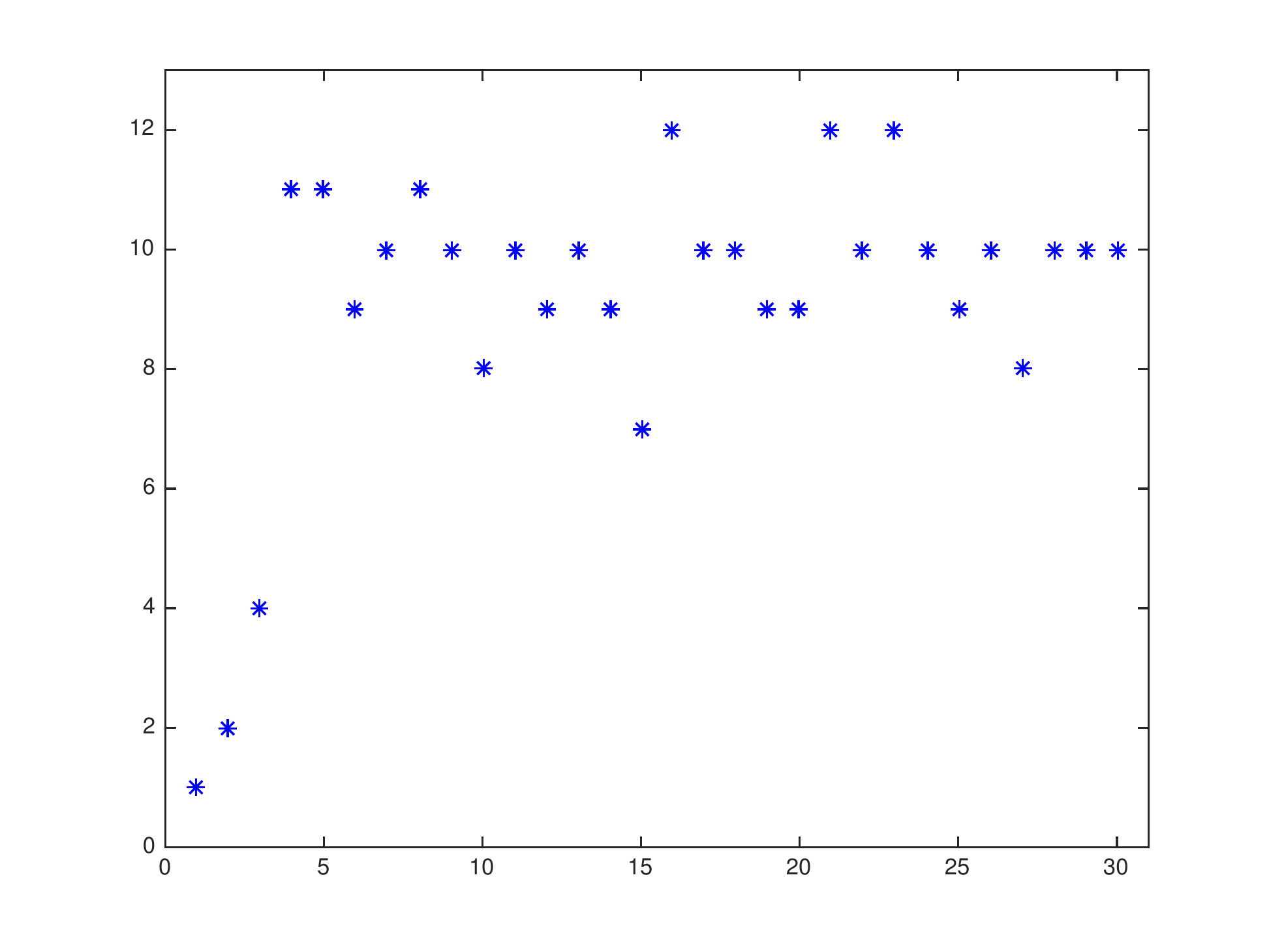}
\caption{Number of launches for the first calendar}
\label{launches1}
\end{center}
\end{figure}

Consider the calendar which number of launches per year is represented on figure \ref{launches1}. They are often 10 launches a year. The ASA algorithm ran for 500 iterations. It took approximatively 8 hours. The results are reported in table \ref{first}.
\begin{table}[tp]
\centering
\begin{tabular}{l l l}\hline
SRM storage capacity & 4 & 8\\
\hline
ASA & 2,158,990 & 1,938,512 \\
Na\"{i}ve policy & 2,336,100 & 2,530,000\\
\hline
\end{tabular}
\label{first}
\caption{Comparison between the performances of the op\-ti\-mal policy computed by the ASA and the na\"{i}ve policy for the first calendar}
\end{table}
The ASA returns better than na\"{i}ve policies. Moreover, storing 8 SRMs leads to a better result (the gain is 23.38\% compared to the na\"{i}ve policy, and it is 10.21\% compared to the capacity of 4 SRMs). Again, the optimal policy computed by ASA does not have a trivial form. It is worth wondering if this conclusion depends on the calendar. 
\subsubsection{Second calendar}
\begin{figure}[tp]
\centering
\includegraphics[scale=0.45]{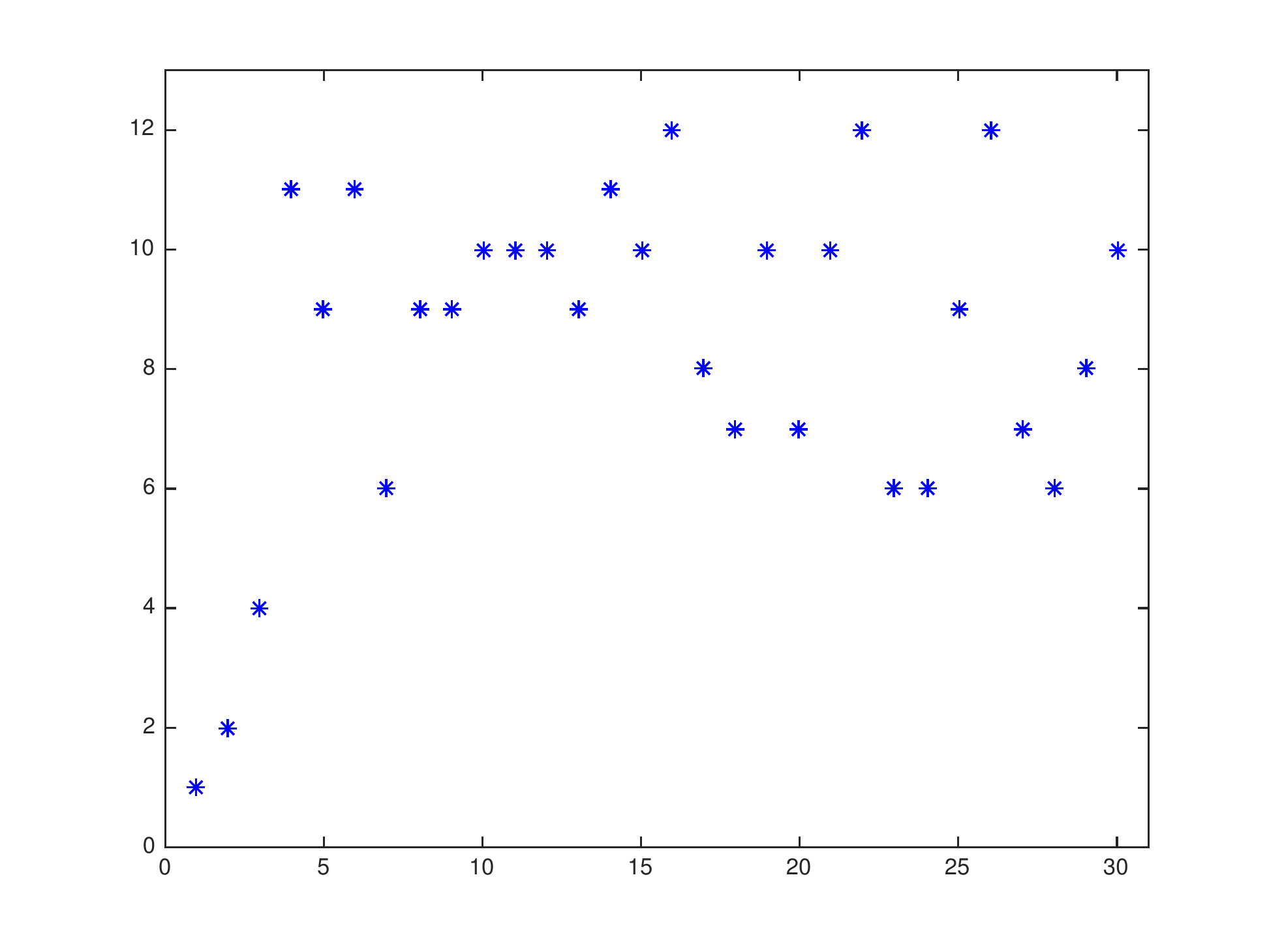}
\caption{Number of launches for the second calendar}
\label{launches2}
\end{figure}

Consider now the calendar on figure \ref{launches2}. The number of launches per year is a little more scattered than in previous example. The results presented in table \ref{second} were obtained with 500 iterations of the ASA.
\begin{table}[tp]
\centering
\begin{tabular}{l l l}\hline
SRM storage capacity & 4 & 8\\
\hline
ASA & 1,452,914 & 1,539,095 \\
Na\"{i}ve policy & 2,452,900 & 2,646,900\\
\hline
\end{tabular}
\label{second}
\caption{Comparison between the performances of the optimal policy computed by the ASA and the na\"{i}ve policy for the second calendar}
\end{table}
In this case, the capacity of 4 SRMs leads to a lower cost (5.6~\% gain). But the gain compared to the na\"{i}ve policy is remarkable: it is 68.83\% (and 71.98\% when the SRM storage capacity is 8 units). The performance depends naturally on the calendar. Actually in this case, there are some years where the number of launches to be performed is comprised between 6 and 8. This leads to lower lateness and storage costs. However, what is more important is that the optimal SRM storage capacity also depends on the calendar. The fact that the number of launches is scattered does not have an influence on this. Indeed, if we take a regular calendar of 10 launches from year 5 to year 30, storing 4 units at most induces a lower cost: the gain is 22.14\% when compared to the 8-unit storage.
\subsubsection{Third calendar}
\begin{figure}[tp]
\centering
\includegraphics[scale=0.45]{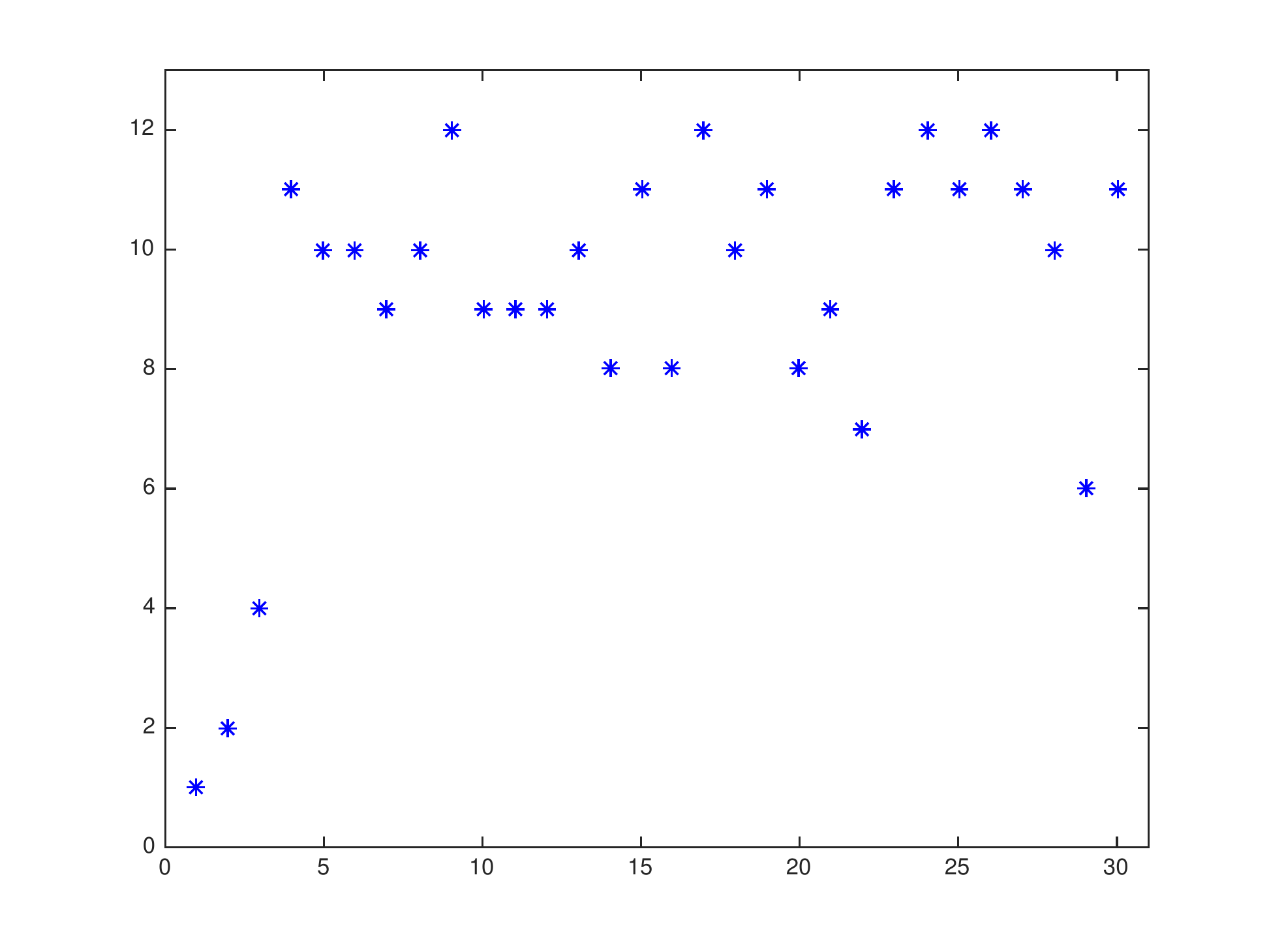}
\caption{Number of launches for the third tested calendar}
\label{launches3}
\end{figure}

For the calendar represented on figure \ref{launches3} and 500 iterations of the ASA, the results obtained are presented in table \ref{third}.
\begin{table}[tp]
\centering
\begin{tabular}{l l l}\hline
SRM storage capacity & 4 & 8\\
\hline
ASA & 2,022,800 & 2,035,200 \\
Na\"{i}ve policy & 2,362,000 & 2,555,700\\
\hline
\end{tabular}
\label{third}
\caption{Comparison between the performances of the optimal policy computed by the ASA and the na\"{i}ve policy for the third calendar}
\end{table}
Storing 4 SRMs at most induces a lower cost. However, the gain compared to an 8-unit storage is only 0.61\% (it is 14.36\% when compared to the na\"{i}ve policy). Taking into account the variance of the costs, one may conclude that, in this case, the two scenarios lead to similar performances.

With these three examples, one sees that the question of the optimal SRM storage capacity is not trivial at all. Indeed, it seems to be impossible to answer this question with only prior knowledge. One may consider static costs to have a more accurate answer to this problem (for instance exploitation costs, a larger warehouse being naturally more expensive). However, it seems that a maximum capacity of 4 SRMs is a better choice in general.
\section{Conclusion}
Using the MDP framework together with a simulation based algorithm, we performed the optimization of the launcher integration process. Several problems had to be addressed, such as different time scales, state space reduction computation speed or numerical representation of numbers.

Given a launch calendar, optimal policies are computed and stored in the form of a matrix. To apply the optimal policy in practice, the controller looks up the best action to select in the matrix, given the current state of the process and the current time. Such policies do not have a trivial form and cannot be easily explained. They lead to up to 70\% gain compared to trivial policies prescribing a manufacturing rate corresponding to the exact number of launches to be performed in the year.

To address real life optimization problem regarding launcher operations, one should now work with a calendar that is known only two years advance. One option that would fit the present framework would be to consider that the calendar of year $n+2$ is randomly drawn at the beginning of year $n$. Further exchange with practitioners is required to derive realistic distributions for such calendars. Another possible extension of interest to Airbus DS is to model the production of subassemblies with more detail, inducing longer delays for selecting the production rates.

\appendix
\section{Instructions for MRAS and ASA algorithms}\label{algos}
Instructions for the MRAS are given on algorithm \ref{MRASa} and for the ASA on algorithm \ref{ASAa}.
They use the following notation.
For some set $A$, $\mathbf{1}_A$ denotes its indicator function, that is $\mathbf{1}_A(x) = 1$ if $x \in A$ and $\mathbf{1}_A(x) = 0$ otherwise.
The upper integer part of $x$ is denoted $\lceil x\rceil$.
Let $\varepsilon > 0$ and define the function $\mathcal{I}$ on $\mathbb{R} \times \mathbb{R}$ by
\begin{equation*}
\mathcal{I}(x,\chi) = 
\left\{ \begin{array}{lcl}
0 &\text{ if }& x \geq \chi + \varepsilon,\\
\dfrac{\chi + \varepsilon - x}{\varepsilon} &\text{ if }& \chi < x < \chi + \varepsilon,\\
1 &\text{ if }& x \leq \chi.
\end{array} \right. 
\end{equation*}
For all states $x_i$, all decisions $a_j$ and all times $t$, let $\Pi_{i,j}(t)$ the set of policies which prescribe action $a_j$ in state $x_i$ at time $t$. For a policy $\pi$, let
\begin{equation*}
f(\pi,P) = \displaystyle \prod_{t = 1}^{30} \prod_{i = 1}^{2281} \prod_{j = 1}^{343} P(i,j,t)^{\mathbf{1}_{\Pi_{i,j}(t)}(\pi)}
\end{equation*}
and, for $\lambda \in (0;1)$, let
\begin{equation*}
\mathbf{f}(\pi,P) = (1 - \lambda)f(\pi,P) + \lambda f(\pi,P_0)
\end{equation*}
with $P_0$ a given probability matrix.

\begin{algorithm}[tp]
\caption{MRAS algorithm}
\begin{algorithmic}
\REQUIRE initial 3-dimensional probability matrix $P_0$, $\rho_0 \in (0;1]$, $\varepsilon > 0$, $N_0 \geq 2$, $M_0 \geq 1$, $\alpha > 1$, $\beta > 1$, $\lambda \in (0;1)$, $\nu \in (0;1]$, initial state $x_0$, iteration count $k = 0$, limit number of iterations $K > 0$.

\WHILE{$k < K$}
	\FOR{$1\leq n\leq N_k$} 
		\STATE -- Draw policy $\pi^n$ from matrix $P_0$ with probability $\lambda$ and matrix $P_k$ with probability $1 - \lambda$.
		\STATE -- Simulate $M_k$ trajectories with $x_0$ as initial state using policy $\pi^n$.
		\STATE -- Compute the cost $V_{k,m}^n$ generated for the trajectories $1\leq m\leq M_k$.
		\STATE -- $\overline{V}_k^n \leftarrow  \dfrac{1}{M_k}\displaystyle\sum_{m=1}^{M_k} V_{k,m}^n$.
	\ENDFOR
	\STATE -- Arrange the $\overline{V}_k^n$ in descending order to get a sequence $(\overline{V}_k^{(n)})$ such that $\overline{V}_k^{(1)} \geq \ldots \geq \overline{V}_k^{(N_k)}$.
	\STATE -- $\gamma_k(\rho_k,N_k) \leftarrow \overline{V}_k^{(\lceil (1-\rho_k)N_k\rceil)}$.
	\IF{$k = 0$ or $\gamma_k(\rho_k,N_k) \leq \overline{\gamma}_{k-1} - \varepsilon$}
		\STATE -- $\overline{\gamma}_k \leftarrow  \gamma_k(\rho_k,N_k)$.
		\STATE -- $\rho_{k+1} \leftarrow \rho_k$.
		\STATE -- $N_{k+1} \leftarrow  N_k$.
		\STATE -- $\boldsymbol{\pi}^*_k \leftarrow $ policy $\pi^n$ such that $\overline{V}^{(n)}_k = \gamma_k(\rho_k,N_k)$.
	\ELSIF{there exists a smallest integer $\nu > \lceil (1-\rho_k)N_k\rceil$ such that $\overline{V}_k^{(\nu)} \leq \overline{\gamma}_{k-1} - \varepsilon$}
		\STATE -- $\overline{\gamma}_k \leftarrow  \overline{V}_k^{(\nu)}$.
		\STATE -- $\rho_{k+1} \leftarrow  1 - \dfrac{\nu}{N_k}$.
		\STATE -- $N_{k+1} \leftarrow  N_k$.
		\STATE -- $\boldsymbol{\pi}^*_k\leftarrow $ policy $\pi^n$ such that $\overline{V}^{(n)}_k = \overline{V}_k^{(\nu)}$.
	\ELSE
		\STATE -- $\overline{\gamma}_k \leftarrow  \overline{V}_k^n$ where $\pi^n = \boldsymbol{\pi}^*_{k-1}$.
		\STATE -- $\rho_{k+1} \leftarrow \rho_k$.
		\STATE -- $N_{k+1} \leftarrow  \lceil \alpha N_k \rceil$.
		\STATE -- $\boldsymbol{\pi}^*_k \leftarrow \boldsymbol{\pi}^*_{k-1}$.
	\ENDIF
	\IF{$\forall n \in \{1,\ldots,N_k\}$, $\dfrac{\mathrm{e}^{-k\overline{V}_k^n}}{\mathbf{f}(\pi^n,P_k)}\mathcal{I}(\overline{V}_k^n,\overline{\gamma}_k) = 0$}
		\STATE -- $P_{k+1} \leftarrow  P_k$.
	\ELSE
		\STATE -- 
$\hat{P}_{k+1}(i,j,t) \leftarrow $ \\\qquad\qquad$\dfrac{\sum_{n=1}^{N_k} \frac{\mathrm{e}^{-k\overline{V}_k^n}}{\mathbf{f}(\pi^n,P_k)}\mathcal{I}(\overline{V}_k^n,\overline{\gamma}_k)\mathbf{1}_{\Pi_{i,j}(t)}(\pi^n)}{\sum_{n=1}^{N_k}{\frac{\mathrm{e}^{-k\overline{V}_k^n}}{\mathbf{f}(\pi^n,P_k)}\mathcal{I}(\overline{V}_k^n,\overline{\gamma}_k)}}.$
		\STATE -- $P_{k+1} \leftarrow  \nu\hat{P}_{k+1} + (1 - \nu)P_k$.
		\STATE -- $M_{k+1} \leftarrow  \lceil \beta M_k \rceil$.
		\STATE --  $k\leftarrow k+1$.
	\ENDIF
\ENDWHILE

\ENSURE matrix $P_K$.
\end{algorithmic}
\label{MRASa}
\end{algorithm}

\begin{algorithm}[tp]
\caption{ASA algorithm}
\begin{algorithmic}
\REQUIRE initial 3-dimensional probability matrix $P_0$, $\alpha_0 = 100^{-0.501}$, $\beta_0 = 1$, temperature parameter $T_0 > 0$, $N_0 > 0$, $M_0 > 0$, initial state $x_0$, iteration count $k = 0$, limit number of iterations $K > 0$.

\WHILE{$k < K$}
	\FOR{$1\leq n\leq N_k$} 
		\STATE -- Draw policy $\pi^n$ from matrix $P_0$ with probability $\beta_k$ and matrix $P_k$ with probability $1 - \beta_k$.
		\STATE -- Simulate $M_k$ trajectories with $x_0$ as initial state using policy $\pi^n$.
		\STATE -- Compute the cost $V_{k,m}^n$ generated for the trajectories $1\leq m\leq M_k$.
		\STATE -- $\overline{V}_k^n \leftarrow  \dfrac{1}{M_k}\displaystyle\sum_{m=1}^{M_k} V_{k,m}^n$.
	\ENDFOR
	\STATE -- 
$\hat{P}_{k+1}(i,j,t) \leftarrow \dfrac{\sum_{n=1}^{N_k} \frac{\exp(-\overline{V}_k^n T_k^{-1})}{\mathbf{f}(\pi^n,P_k)}\mathbf{1}_{\Pi_{i,j}(t)}(\pi^n)}{\sum_{n=1}^{N_k}{\frac{\exp(-\overline{V}_k^n T_k^{-1})}{\mathbf{f}(\pi^n,P_k)}}}.$
	\STATE -- $P_{k+1} \leftarrow \nu\hat{P}_{k+1} + (1 - \nu)P_k$.
	\STATE -- $M_{k+1} \leftarrow \max(M_0,[1.10\log^3(k)])$.
	\STATE -- $N_{k+1} \leftarrow \max(N_0,[k^{0.501}])$.
	\STATE -- $\alpha_{k+1} \leftarrow (k + 100)^{-0.501}$.
	\STATE -- $\beta_{k+1} \leftarrow (\sqrt{k + 1})^{-1}$.
	\STATE -- $T_k \leftarrow T_0(\log(k + \exp(1)))^{-1}$.
	\STATE -- $k\leftarrow k+1$.
\ENDWHILE

\ENSURE matrix $P_K$.
\end{algorithmic}
\label{ASAa}
\end{algorithm}

\end{document}